\RequirePackage{ifpdf}

\RequirePackage{fixltx2e}
\RequirePackage{fix-cm}

\documentclass[%
floatfix,
showkeys,
nofootinbib, %
superscriptaddress, %
]{revtex4-1}

\usepackage{iftex}
\ifpdftex
\usepackage[noTeX]{mmap}	
\fi

\usepackage{textcomp}

\usepackage[T1,T2A]{fontenc}

\usepackage[utf8]{inputenx}
\input{ix-utf8enc.dfu}

\usepackage[german,russian,english]{babel}

\usepackage{amsmath}
\usepackage{amssymb}

\usepackage{mathtools}
\mathtoolsset{
showonlyrefs,
mathic = true
}

\allowdisplaybreaks

\usepackage{hyperref}
\hypersetup{backref,
 colorlinks=false}
\hypersetup{pdfborder=0 0 0}

\usepackage{microtype}
\UseMicrotypeSet[protrusion]{alltext}

\usepackage{graphicx}

\usepackage[scanall]{psfrag}

\usepackage{listings}
\usepackage{listingsutf8}
\usepackage[frozencache]{minted}

\usepackage{fixltx2e}

\usepackage{nicefrac}

\makeatletter
\def\ps@pprintTitle{%
     \let\@oddhead\@empty
     \let\@evenhead\@empty
     \let\@oddfoot\@empty
     \let\@evenfoot\@oddfoot}
\makeatother

\usepackage{booktabs}

\begin{document}

\graphicspath{{image/linear-phase/en/}{image/linear-phase/}{image/}}

\title{Numerical integrating of highly oscillating functions: effective stable algorithms in case of linear phase}

\author{Leonid A. Sevastianov}
\email{sevastianov-la@rudn.ru}
\affiliation{
  Peoples' Friendship University of Russia (RUDN University),\\
  6 Miklukho-Maklaya St, Moscow, 117198, Russian Federation}
\affiliation{
  Joint Institute for Nuclear Research\\
  6 Joliot-Curie, Dubna, Moscow region, 141980, Russian Federation}

\author{Konstantin P. Lovetskiy}
\email{lovetskiy-kp@rudn.ru}
\affiliation{
  Peoples' Friendship University of Russia (RUDN University),\\
  6 Miklukho-Maklaya St, Moscow, 117198, Russian Federation}

\author{Dmitry S. Kulyabov}
\email{kulyabov-ds@rudn.ru}
\affiliation{
  Peoples' Friendship University of Russia (RUDN University),\\
  6 Miklukho-Maklaya St, Moscow, 117198, Russian Federation}
\affiliation{
  Joint Institute for Nuclear Research\\
  6 Joliot-Curie, Dubna, Moscow region, 141980, Russian Federation}

\begin{abstract}
  A practical and simple stable method for calculating Fourier
  integrals is proposed, effective both at low and at high
  frequencies. An approach based on the fruitful idea of Levin, to use
  of the collocation method to approximate the slowly oscillating part
  of the antiderivative of the desired integral, allows reducing the
  calculation of the integral of a highly oscillating function (with
  a linear phase) to solving a system of linear algebraic equations
  with a three-diagonal triangular or five-diagonal band Hermitian
  matrix.
  The choice of Gauss-Lobatto grid nodes as collocation points makes
  it possible to use the properties of discrete ``orthogonality'' of
  Chebyshev differentiation matrices in physical and spectral
  spaces. This is realised in increasing the efficiency of the
  numerical algorithm for solving the problem. The system
  pre-conditioning procedure leads to significantly less cumbersome
  and more economical calculation formulas.  To avoid possible
  numerical instability of the algorithm, we proceed to the solution
  of a normal system of linear algebraic equations.
\end{abstract}


\maketitle

\section{Introduction}

The initial formulation of the method of numerical integration of
highly oscillating functions by Levin and his followers suggests a
possible ambiguity in finding the antiderivative: any solution to the
differential equation without boundary (initial) conditions can be used
to calculate the desired value of the integral.

Levin's approach \cite{1} to the integration of highly oscillating
functions consists in the transition to the calculation of the
antiderivative function from the integrand using the collocation
procedure in physical space. In this case, the elements of the
degenerate \cite{2} differentiation matrix of the collocation method
\cite{3} are a function of the coordinates of the grid points, the matrix
elements are calculated using very simple formulas. In books \cite{3,4}
various options for the implementation of this method are considered,
many applied problems are solved.

The method proposed by Levin both in the one-dimensional and in the
multidimensional case was published by him in articles \cite{1,5}, and
then he was thoroughly studied in \cite{6}. The method is presented in
great detail in the famous monograph \cite{4}, which describes the
evolution of numerical methods for integrating highly oscillating
functions over the past fifteen years.

There are a large number of works using various approaches in order to
propose fast and stable methods for solving systems of linear algebraic
equations (SLAE) that arise when implementing the collocation method.
However, many of them \cite{7,8,9} encounter difficulties in solving the
corresponding systems of linear equations.

In particular, the use in specific implementations of the Levin
collocation method in the physical space of degenerate Chebyshev
differentiation matrices, which also have eigenvalues differing by
orders of magnitude, makes it impossible to construct a stable numerical
algorithm for solving the resulting SLAEs. The approach to solving the
differential equation of the Levin method, described in \cite{8,10,11},
is based on the approximation of the solution, as well as the integrand
phase and amplitude functions in the form of expansion into finite
series in Chebyshev polynomials. Moreover, to improve the properties of
the algorithms, and hence the matrices of the corresponding SLAEs,
three-term recurrence relations are used that connect the values of
Chebyshev polynomials of close orders. However, these improvements are
not enough to ensure stable calculation of integrals with large matrix
dimensions.

In our work, we consider a method of constructing a primitive, based on
the spectral representation of the desired function.

We propose increasing the efficiency of the algorithm by reducing the
corresponding system of linear equations to a form that is always
successfully solved using the LU-decomposition method with partial
selection of the leading element.

Consider the integral that often occurs in Fourier analysis--in
applications related to signal processing, digital images, cryptography
and many other areas of science and technology

\begin{equation}
  \label{eq:1}
  I_{\omega}(f) =
  \int_{a}^{b}f(x)e^{i\omega g(x)} dx.
\end{equation}

In accordance with the Levin method, the calculation of this integral
reduces to solving an ordinary differential equation
\begin{equation}
  \label{eq:2}
  p'( x ) + i\omega g' \left( x \right)p\left( x \right) = f\left( x \right),\ x \in \lbrack a,b\rbrack.
\end{equation}

As argued in \cite{1}, the system~(\ref{eq:2}) has a particular
solution which is not highly oscillatory, and we shall look for an
approximation to this particular solution by collocation with `nice'
functions, e.g.  polynomials. If the unknown function
\(p\left( x \right)\) is a solution of Eq.~(\ref{eq:2}), then the
result of integration can be obtained according to the formula

\begin{equation}
  \label{eq:3}
  I_{\omega}\left( f,g \right) =
  \int_{a}^{b}{\left( p' \left( x \right) + i\omega g' \left( x
      \right)p\left( x \right) \right)
    e^{{i\omega g}\left( x \right)}{dx}} =
  p\left( b \right)e^{{i\omega g}\left( b \right)} - p\left( a \right)e^{{i\omega g}\left( a \right)}.
\end{equation}

Below we will consider the special case of integration of a highly
oscillating function with a linear phase, reduced to the standard form.

\begin{equation}
  \label{eq:4}
  I_{\omega}\left\lbrack f \right\rbrack
  = \int_{-1}^{1} f(x)e^{{i\omega x}}{dx}
  = p (1) e^{{i\omega}}
  - p\left( - 1 \right)e^{- i\omega}.
\end{equation}

This can be justified, in particular, by the fact that in many
well-known publications~\cite{9,12,13} stable transformations are discussed
in detail, which make it possible to proceed from a general integral
with a nonlinear phase to an integral in standard form (on the
interval $[-1 , 1]$) with a linear phase.

In the paper by Levin~\cite{1}, to automatically exclude
the highly oscillating component \(ce^{- i\omega x}\) of the
general solution
\({\ p}\left( x \right) = p_{0}\left( x \right) + ce^{- i\omega g(x)}\),
it is proposed to search for a numerical solution~(\ref{eq:2}) based on the
collocation method, using its expansion in a basis of slowly oscillating
functions, rather than using difference schemes (or methods of the
Runge--Kutta type).

In this case, the following statement is true~\cite{2}:

\textbf{Statement.} The solution of Eq.~(\ref{eq:2}) obtained using the Levin
collocation method is a slowly oscillating function
\(\mathcal{O} (\omega^{- 1})\) for \(\omega \gg 1\).

\section{Approximation of the antiderivative. Calculation method}

Let us consider in more detail the problem of finding the
antiderivative integrand, or rather, the approximating polynomial
\(p(x)\), satisfying condition~(\ref{eq:2}) in a given number of
points on the interval $[-1,1]$. Consider the spectral method of
finding an approximating function in the form of expansion in a finite
series

\begin{equation}
  \label{eq:5}
  p\left( x \right) = \sum_{k = 0}^{n}{c_{k}T_{k}(x)},
  \quad
  x \in \lbrack - 1,1\rbrack
\end{equation}
in the basis of Chebyshev polynomials of the first kind
\(\left\{ T_{k}(x) \right\}_{k = 0}^{\infty}\), defined in the Hilbert
space of functions on the interval $[-1,1]$.

The application of the collocation method to solve the problem
\(p' (x) + i\omega p(x) = f(x)\) leads to the need to fulfill the
following equalities for the desired
coefficients\(\ c_{k},\ k = 1,\ldots,n\)

\begin{equation}
  \label{eq:6}
  \sum_{k = 0}^{n}{c_{k}T_{k}' \left( x_{j} \right)} + i\omega\sum_{k = 0}^{n}{c_{k}T_{k}\left( x_{j} \right)} = f\left( x_{j} \right),\ j = 0,\ldots,n
\end{equation}
at the collocation points \(\left\{ x_{0},x_{1},\ldots,x_{n} \right\}\).

The last statement is equivalent to the fact that the coefficients
\(c_{k},\ k = 0,\ldots,n\) should be a solution to the system of linear
algebraic equations of the collocation method:
\begin{equation}
  \label{eq:7}
  \left\{
    \begin{gathered}
      p' \left( x_{0} \right) + i\omega p(x_{0}) = f\left( x_{0} \right), \\
      p' \left( x_{1} \right) + i\omega p(x_{1}) = f\left( x_{1} \right), \\
      \ldots \\
      p' \left( x_{n} \right) + i\omega p(x_{n}) = f\left( x_{n} \right). \\
    \end{gathered}
  \right.
\end{equation}

We represent the values of the derivative of the desired function
(polynomial) at the collocation points in the form of the product
\({Dp = p'}\) of the matrix \({D}\) by the vector of
values of \({p}\). Recall that the Chebyshev differentiation
matrix \({D}\) has the standard representation in the physical
space~\cite{3}
\begin{equation}
  \label{eq:8}
  {D}_{{{kj}}}=
  \left\{
    \begin{aligned}
      &\frac{r_{k}}{r_{j}} \left( - 1 \right)^{k + j}/(x_{k} - x_{j}), & k,j = 0,\ldots n,k \neq j \\
      &- \sum_{l = 0,l \neq k}^{n}D_{{kl}}, & k = j.
    \end{aligned}
  \right.
\end{equation}
where
$$
r_{j} =
\left\{
  \begin{matrix}
    2, & j = 0,n \\
    1, & 1,\ldots,n - 1.
  \end{matrix}
\right.
$$

Substituting \({p}^{{'}}{= Dp}\) into Eq.~(\ref{eq:7}) we
reduce it to a system of linear algebraic equations
\begin{equation}
  \label{eq:9}
  \left( {D +}{i\omega }{E} \right){p = f}.
\end{equation}

Here $E$ is an identity matrix, \(f\) is a vector of
values of the amplitude function on the grid. Denote by \({B}\)
the differentiation matrix in the frequency (spectral) space~\cite{14},
whose coefficients are explicitly expressed as
\begin{equation}
  \label{eq:10}
  B_{{{ij}}}{=}\left\{
    \begin{matrix}
      (1/r_{j}) 2 j & \text{if } j > i, i + j \text{ odd}, \\
      0 & \text{ otherwise}
    \end{matrix}
\right.
\end{equation}
where
$$
0 \leq i,j \leq n  \wedge  r_{i} =
\left\{ \begin{matrix}
2 & i = 0 \\
1 & i > 0. \\
\end{matrix}
\right.
$$

Denote by \({T}\) the Chebyshev matrix of mapping a point (vector)
from the space of coefficients to the space of values of the function
\cite{15}. Given that \({p =}{{Tc}}\) is the vector of values of the
desired function (also in physical space), the components of the
derivative vector can be written as \({D}{p =}{{TBc}}\)~\cite{15}. As
a result, we obtain the system of linear algebraic equations
equivalent to system~(\ref{eq:9}),
\begin{equation}
  \label{eq:11}
  \left( {{TBc}} + i\omega{{Tc}} \right)= f
\end{equation}
which is valid for an arbitrary grid on the interval $[-1,1]$. We
write equation~(\ref{eq:11}) in detail
\begin{equation}
  \label{eq:12}
  \begin{bmatrix}
    T_{00} & T_{10} & T_{20} & \vdots & T_{n0} \\
    T_{01} & T_{11} & T_{21} & \vdots & T_{n1} \\
    T_{02} & T_{12} & T_{22} & \vdots & T_{n2} \\
    \ldots & \ldots & \ldots & \ddots & \ldots \\
    T_{0n} & T_{1n} & T_{2n} & \vdots & T_{nn} \\
  \end{bmatrix}
  \left(
    \begin{bmatrix}
      0 & 1 & 0 & 3 & \vdots \\
      & 0 & 4 & 0 & \vdots \\
      &   & 0 & 6 & \vdots \\
      &   &   & \ddots & \vdots \\
      &   &   &   & 0
    \end{bmatrix}
    + {i\omega }{E}
  \right)
  \begin{bmatrix}
    c_{0} \\
    c_{1} \\
    c_{2} \\
    \ldots \\
    c_{n} \\
  \end{bmatrix}
  =
  \begin{bmatrix}
    f_{0} \\
    f_{1} \\
    f_{2} \\
    \ldots \\
    f_{n}
  \end{bmatrix}
\end{equation}
where to reduce the formulas we used the notation
\(T_{{kj}} = T_{k}\left( x_{j} \right),\ k,j = 0,\ldots,n\).

The product of a non-degenerate matrix \(T\) by a non-degenerate
triangular matrix \({B +}i\omega{E}\) is a non-degenerate
matrix. Therefore, the system of linear algebraic equations~(\ref{eq:12}) has a
unique solution.

\paragraph{Statement 1}
The solution of this system of linear algebraic equations with respect
to the coefficients allows us to approximate the antiderivative function
in the form of a series~(\ref{eq:5}) and calculate the approximate value of the
integral by formula~(\ref{eq:4}).

\section{Modification of the calculation method}

System~(\ref{eq:12}) is valid for an arbitrary grid on the interval $[-1,1]$.
However, consideration of the collocation problem on a Gauss--Lobatto
grid allows significant simplification of this system of linear
algebraic equations. First, we multiply the first and last equations
from~(\ref{eq:12}) by \(1/\sqrt{2}\) to obtain an equivalent ``modified'' system
with a new matrix \(\tilde{{T}}\) (instead of
\({T}\)), which is good because it has the property of discrete
``orthogonality'' and, therefore, is non-degenerate. Therefore,
multiplying it on the left by its transposed one gives the diagonal
matrix:
$$
{\tilde{{T}}}^{{T}} \tilde{{T}}{=}\begin{bmatrix}
n & 0 & 0 & \vdots & 0 \\
0 & n/2 & 0 & \vdots & 0 \\
0 & 0 & n/2 & \vdots & 0 \\
\ldots & \ldots & \ldots & \ddots & \ldots \\
0 & 0 & 0 & \vdots & n \\
\end{bmatrix}.
$$

We use this property and multiply the reduced (modified) system~(\ref{eq:12}) on
the left by the transposed matrix
\({\tilde{{T}}}^{{T}}\), thereby reducing it to the
upper triangular form. Indeed, the matrix of the resulting system is
calculated as the product of the diagonal matrix by the triangular
matrix, which, in turn, is the sum of the Chebyshev differentiation
matrix in the spectral space and the diagonal matrix.

Since the matrix \({\tilde{{T}}}^{{T}}\) is
non-degenerate, the new system of linear algebraic equations is
equivalent to system~(\ref{eq:12}) and has a unique solution.

Taking into account the specific values of the Chebyshev polynomials
\emph{on the Gauss-Lobatto grid}~\cite{16}, simplifies the system,
bringing it to the form
\begin{equation}
  \label{eq:13}
{{Ac}} = \begin{bmatrix}
{i\omega } & 1 & 0 & 3 & \vdots & n - 1 \\
0 & {i\omega } & 2 & 0 & \vdots & 0 \\
0 & 0 & {i\omega } & 3 & \vdots & n - 1 \\
0 & 0 & 0 & {i\omega } & \vdots & 0 \\
\ldots & \ldots & \ldots & \ldots & \ddots & n - 1 \\
0 & 0 & 0 & 0 & \vdots & {i\omega } \\
\end{bmatrix}\begin{bmatrix}
c_{0} \\
c_{1} \\
c_{2} \\
c_{3} \\
\ldots \\
c_{n} \\
\end{bmatrix} = \begin{bmatrix}
{\tilde{f}}_{0}/2 \\
{\tilde{f}}_{1} \\
{\tilde{f}}_{2} \\
{\tilde{f}}_{3} \\
\ldots \\
{\tilde{f}}_{n}/2 \\
\end{bmatrix}
\end{equation}
where
\({\tilde{f}}_{j} = \frac{1}{n}\sum_{k = 0,n}^{''}{T_{j}\left( x_{k} \right)f\left( x_{k} \right)},\ j = 0,\ldots,n\)
and symbol \(\Sigma''\) denotes a sum in which the first and last
terms are additionally multiplied by 1/2.

By the Kronecker--Capelli theorem, the system of linear algebraic
equations~(\ref{eq:13}) with a square matrix and a non-zero determinant is not
only solvable for any vector of the right-hand side, but also has a
unique solution.

\paragraph{Statement 2} For \(\left| \omega \right| > n\) the SLAE~(\ref{eq:13})
has a stable solution.

\emph{The reverse course of the Gauss method for solving system~(\ref{eq:13}) can
lead to accumulation of errors for $n >|\omega|$.}

\paragraph{Statement 3} To solve system~(\ref{eq:13}), no more than
(\(\sim n^{2}/4\)) operations of addition/subtraction and
multiplication/division with a floating point are required.

\section{Efficient method for solving the problem}

To increase the efficiency of the method for solving the system of
linear algebraic equations~(\ref{eq:13}), we simplify it by reducing to a
triangular band three-diagonal form. To do this, we multiply system~(\ref{eq:13})
on the left by a band non-degenerate matrix with unit diagonal elements
and equal to $-1$ elements on the second upper codiagonal
\begin{equation}
  \label{eq:14}
  Q_{ij}=
  \left\{
    \begin{matrix}
      1 & \text{if } i = j, \\
      -1 & j = i + 2.
    \end{matrix}
  \right.
\end{equation}

As a result, the matrix of system~(\ref{eq:13}) takes on a fairly simple form of
a band supra-diagonal matrix with non-zero elements only on the main
diagonal and on two upper codiagonals:
\begin{equation}
  \label{eq:15}
  QAc = Gc =
  \begin{bmatrix}
    i \omega & 1 & - i \omega & 0 & 0 & 0 & 0 \\
    0 & i \omega & 2 & - i \omega & 0 & 0 & 0 \\
    0 & 0 & i \omega & 3 & - i \omega & 0 & 0 \\
    \ldots & \ldots & \ldots & \ddots & \ldots & - i \omega & \ldots \\
    0 & 0 & 0 & 0 & i \omega & n - 2 & - i \omega \\
    0 & 0 & 0 & 0 & 0 & i \omega & n - 1 \\
    0 & 0 & 0 & 0 & 0 & 0 & i \omega \\
  \end{bmatrix}\begin{bmatrix}
    c_{0} \\
    c_{1} \\
    c_{2} \\
    \ldots \\
    c_{n - 2} \\
    c_{n - 1} \\
    c_{n} \\
  \end{bmatrix}{= Q}\begin{bmatrix}
    {\tilde{f}}_{0}/2 \\
    {\tilde{f}}_{1} \\
    {\tilde{f}}_{2} \\
    \ldots \\
    {\tilde{f}}_{n - 2} \\
    {\tilde{f}}_{n - 1} \\
    {\tilde{f}}_{n}/2 \\
  \end{bmatrix}
\end{equation}

Since the matrix \({Q}\) is non-degenerate, the system of linear
algebraic equations~(\ref{eq:15}) is equivalent to system~(\ref{eq:13}), therefore, it has
a unique solution.

Algorithms for solving systems of linear equations such as the Gauss
method or the LU-decomposition work well when the matrix of the system
has the property of diagonal dominance. Otherwise, standard solution
methods lead to the accumulation of rounding errors. A stable solution
to the system is provided by the LU-decomposition method with a partial
choice of a leading element.

A solution to system~(\ref{eq:17}) can still be unstable for the same reason as
in the case of system of equations~(\ref{eq:13}).

Passing to the solution of the normal system \cite{17}, that is, to the
problem of minimizing the residual
\(\left\| {{Ac}}{-}\tilde{{f}} \right\|^{2}\),
multiplying the system of equations~(\ref{eq:15}) on the left by the Hermitian
conjugate matrix
\begin{equation}
  \label{eq:16}
{A}^{\dagger}{{Ac}} = {A}^{\dagger}\tilde{{f}}
\end{equation}
we transform the matrix of the system~(\ref{eq:15}) to the five-diagonal form:
\begin{equation*}
  {A}^{\dagger}A =
  \begin{bmatrix}
    -\omega^2 & -i\omega & \omega^2 & 0  & \vdots & 0 & 0 \\
    i\omega & 1 -\omega^2 & -3 i \omega & \omega^2 & \vdots & 0 & 0 \\
    \omega^2 & 3 i \omega & 4 - 2 \omega^2 & -5 i \omega  & \vdots & 0 & 0 \\
    0 & \omega^2 & 5 i \omega & 9 - 2\omega^2  & \vdots & 0 & 0 \\
    0 & 0 & \omega^2 & 7 i \omega  & \vdots & 0 & 0 \\
    \cdots & \cdots & \cdots & \cdots  & \ddots & \cdots & \cdots \\
    0 & 0 & 0 & 0 & \ldots & (n - 1)^{2} - 2\omega^2 & - i \omega(2n - 1) \\
    0 & 0 & 0 & 0 & \vdots & i\omega (2n-1) & n^{2} - 2\omega^2 \\
  \end{bmatrix}
\end{equation*}

Although the system of linear equations
\({A}^{\dagger}{{Ac}} = {A}^{\dagger}\tilde{{f}}\) became more filled, since
instead of upper triangle three-diagonal matrix a system of linear
equations with five-diagonal matrix appeared, its computational
properties are cardinally improved. The resulting matrix of a system
of linear algebraic equations is Hermitian, its eigenvalues are real,
and the eigenvectors form an orthonormal system.  The method of
LU-decomposition with a partial choice of the leading element, due to
the properties of the resulting matrix, provides~\cite{17} the
stability of the numerical algorithm for finding the only solution to
the system.

\paragraph{Statement 4} To solve a system of linear algebraic equations
with a band 5-diagonal matrix, the number of required floating-point
operations is of the order of $(19n-29)$ (i.e. \(\mathcal{O}(n)\))~\cite{18}.

\section{Description of the algorithm}

Let us describe the sequence of operations of the presented algorithm
for calculating the integral of a highly oscillating function of the
form~(\ref{eq:1}) with a linear phase.

\paragraph{Input data preprocessing}

\begin{enumerate}
\item
  If the integral is given on the interval \(\lbrack a,b\rbrack\), we
  pass to the standard domain of integration \(\lbrack - 1,1\rbrack\) by
  changing the variables
  \(x = \frac{b - a}{2}t + \frac{b + a}{2},\ t \in \lbrack - 1,1\rbrack\).
\item
  Fill by columns the Chebyshev transformation matrix~(\ref{eq:12}) using only
  one pass of the recursive method for calculating the values of
  Chebyshev polynomials of the first kind of the n-th order.
\end{enumerate}

\paragraph{Antiderivative algorithm}

\begin{enumerate}
\setcounter{enumi}{3}
\item
  Calculate the vector of the right-hand side of system~(\ref{eq:15})
\item
  Fill in the elements of the sparse matrix~(\ref{eq:16}), which depend only on
  the dimension $n$ and the phase value $\omega$.
\item
  If \(\left| {i\omega } \right| > n\), then go to step 6. Otherwise go
  to step 7.
\item
  The matrix of system~(\ref{eq:16}) is a matrix with a diagonal dominance and
  can be stably solved. The solution values at the boundary points are
  used to determine the desired antiderivative values. Go to step 8.
\item
  Multiply relation~(\ref{eq:15}) on the left by the conjugate matrix to obtain a
  Hermitian matrix with diagonal dominance. In this case, to determine
  the values of the antiderivative at the boundary points the normal
  solution is stably determined using the LU-decomposition with a
  partial choice of the leading element.
\item
  We calculate the values of the antiderivative at the ends of the
  interval using the formulas
  \(p\left( 1 \right) = \sum_{j = 0}^{n}c_{j},\ \) and
  \(p\left( - 1 \right) = \sum_{j = 0, j \text{ is even}}^{n}c_{j} -
  \sum_{j = 0, j \text{ is odd}}^{n}c_{j}\).
  The desired value of the integral is obtained using the formula
  \(I(f,\omega) = p\left( 1 \right)e^{{i\omega }} - p\left( - 1 \right)e^{- i\omega}\).
\end{enumerate}

\section{Numerical examples}

\subsection{Example 1}

We give an example of calculating the integral when, for a good
polynomial approximation of a slowly oscillating factor of the
integrand, it is necessary to use polynomials of high degrees.

\begin{equation}
  \label{eq:17}
  I_{\omega}\left\lbrack \frac{1}{x + 2} \right\rbrack =
  \int_{-1}^{1}{\frac{1}{x + 2}e^{{i\omega x}}{dx}}.
\end{equation}

This integral is given by Olver~\cite[p. 6]{19} as an example of the
fact that the GMRES method allows one to calculate the integral much
more accurately than the Levin collocation method. However, in his
article, solving the resulting system of linear algebraic equations
requires \(\mathcal{O(}n^{3})\) operations, as in the Levin collocation
method using the Gaussian elimination algorithm.

The following table shows the values of the integral calculated by us
for various values of the parameter $\omega$ with an accuracy of 17 significant
digits.

\begin{table}
  \caption{The following table shows the values of the integral
    calculated for various values of the parameter $\omega$ with an
    accuracy of 17 significant digits}
  \label{tab:1}
  \centering
\begin{tabular}{lll}
\toprule
$\omega$  & Real part of (\ref{eq:17}) & Imaginary part of (\ref{eq:17})\\
\midrule
$\omega =1$ & 0.9113301035062809891 & -0.1775799622517861791\tabularnewline
$\omega =10$ & -0.07854759997855625023 & -0.04871911238563061052\tabularnewline
$\omega =50$ & -0.00665013790168713 & 0.0129677770647216\tabularnewline
$\omega =100$ & -0.00667389328931381 & 0.00580336592710437\tabularnewline
\bottomrule
\end{tabular}
\end{table}

\begin{figure}
  \centering
  \includegraphics[width=0.8\linewidth]{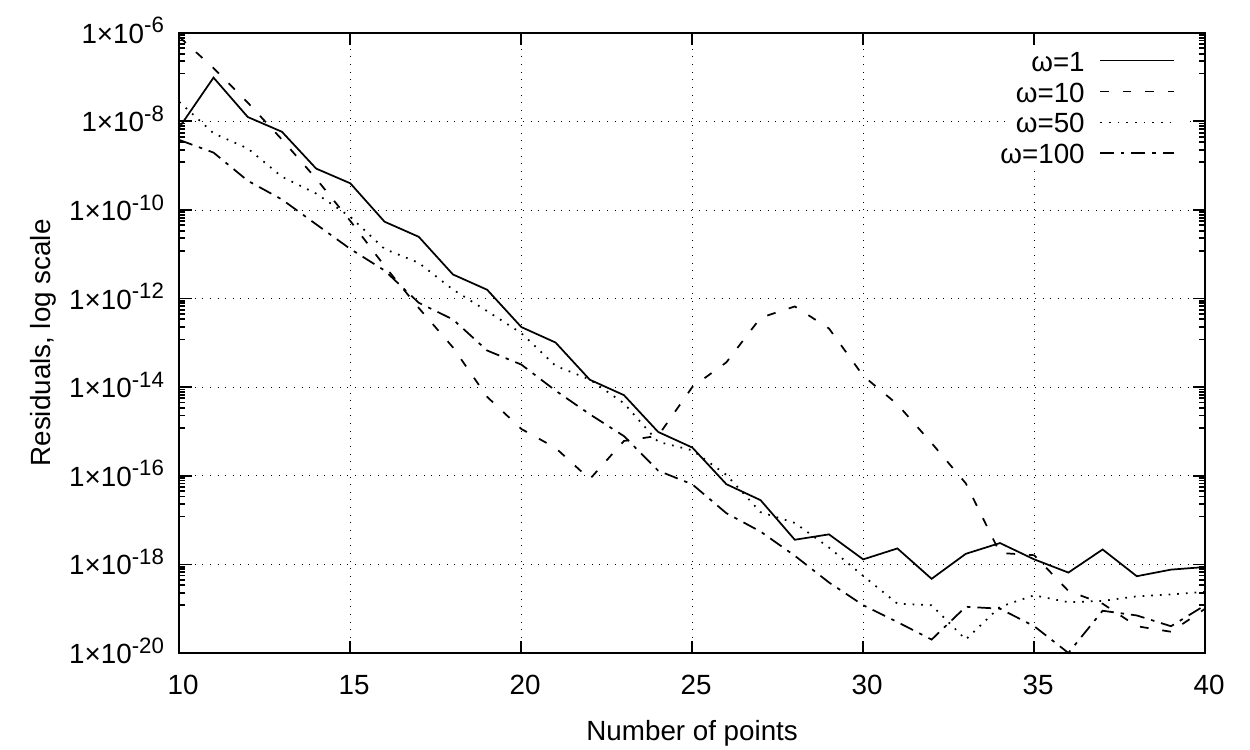}
  \caption{The error in approximating integral~(\ref{eq:17}) for
    different choices of \(\omega\)}
  \label{fig:1}
\end{figure}

A comparison of our results at 40 interpolation points with the results
of~\cite{19} shows a significant gain in accuracy: the deviation from the
exact solution is of the order of $10^{-17}$ compared with
the deviation of the order of $10^{-7}$ in Olver's article.
The proposed algorithm to achieve an accuracy of $10^{-13}$
in the calculation of the integral uses no more than 30 points
(\(n \leq 30)\) for $\omega =1,\ldots, 100$. Moreover, to solve a system of
linear algebraic equations with a band 5-diagonal matrix, the order of
$(19n-29)$ \(\mathcal{O}(n)\) operations is required~\cite{18}.

\subsection{Example 2}

As a second example, we consider the integral
\begin{equation}
  \label{eq:18}
  \int_{- 1}^{1}{\frac{1}{x^{2} + 1}e^{i\omega \sin(x + 1/4)}{dx}}
\end{equation}
from~\cite{19}, where the results of calculating the integrals depending
on the number of approximation points are illustrated (see Fig.~\ref{fig:1}).

To reduce this integral to the (standard form of the Fourier integral)
form of integral with the linear phase, we change the variables
\(y = \sin(x + \frac{1}{4})\).
Then\(\ dx = \frac{1}{\sqrt{1 - y^{2}}}{dy}\),
\(x = \arcsin\left( y \right) - 1/4\), the integration limits are
changed to
[\(- \sin\left( \frac{3}{4} \right), \sin(\frac{5}{4})\)]
and the integral can be written as:
\begin{equation}
  \label{eq:19}
  \int_{- \sin(\frac{3}{4})}^{\sin(\frac{5}{4})}{\frac{1}{\sqrt{1 - y^{2}}{((\arcsin\left( y \right) - \frac{1}{4})}^{2} + 1)}e^{{i\omega y}}}{dy}.
\end{equation}

Let us consider the calculation of this integral for various values of
the parameter $\omega$ using an algorithm that takes into account the linearity
of the phase function. The following table shows the values of the
integral (calculated with high accuracy using exterior programs on the
website \url{https://www.wolframalpha.com/}) for various values of the
parameter $\omega$. The table compares the ``exact'' integration results with
the values calculated using our program on a grid of $90$ points for
various values of the frequency $\omega$. The figure below shows plots
comparing the results compared with the results of~\cite{19},
demonstrating the ``advantage'' in the accuracy of our numerical results
on grids of the same size.

\begin{table}
  \caption{The table shows the values of the integral for various values of the parameter $\omega$}
  \label{tab:2}
  \centering
\begin{tabular}{lll}
\toprule
$\omega$  & Real part of (\ref{eq:19}) & Imaginary part of (\ref{eq:19})\\
\midrule
$\omega =0.1$ & 1.5687504317409 & 0.0337582105322438\tabularnewline
$\omega =1$ & 1.3745907842843 & 0.305184104407599\tabularnewline
$\omega =3$ & 0.311077689499021 & 0.339612459676631\tabularnewline
$\omega =10$ & 0.00266714972608754 & 0.180595659138141\tabularnewline
$\omega =30$ & 0.00706973992290492 & 0.0455774930833239\tabularnewline
$\omega =50$ & -0.00620005944852318 & 0.0155933115982172\tabularnewline
$\omega =100$ & 0.00460104072965418 & -0.00790563176002816\tabularnewline
\bottomrule
\end{tabular}
\end{table}

\begin{figure}
  \centering
  \includegraphics[width=0.8\linewidth]{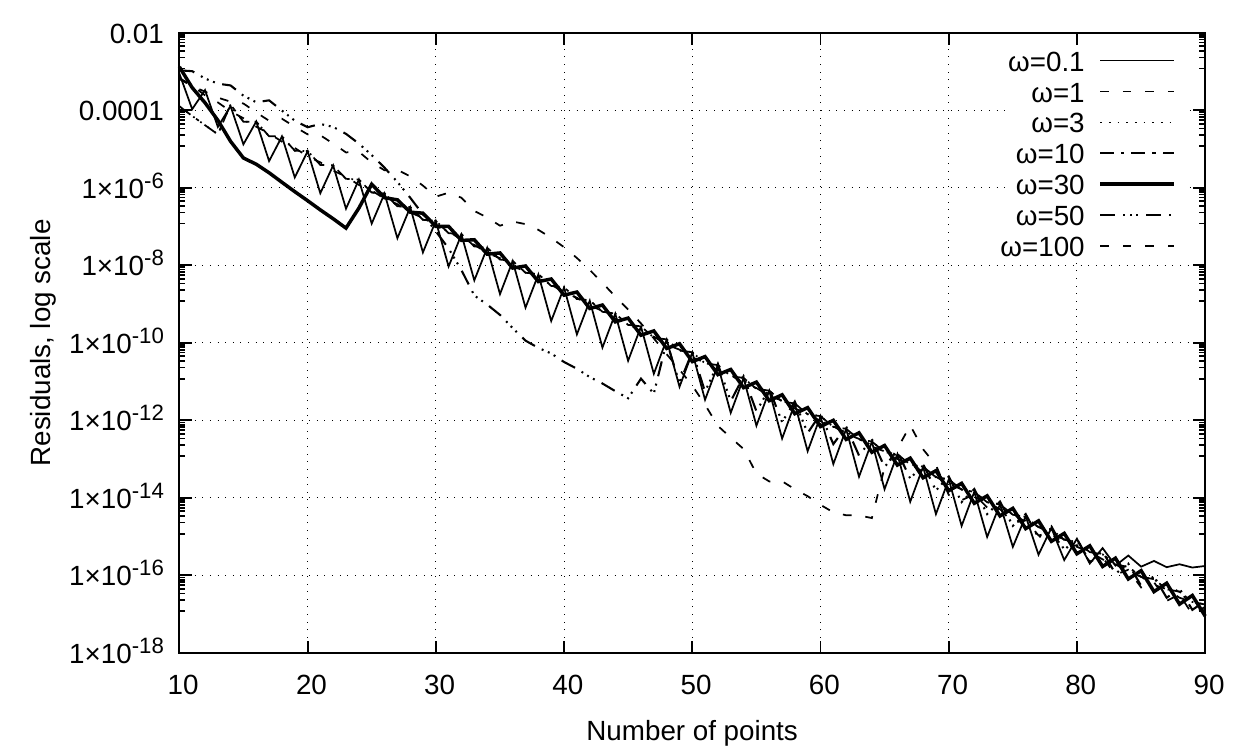}
  \caption{The error in approximating integral~(\ref{eq:19}) for
    different choices of $\omega$}
  \label{fig:2}
\end{figure}

The proposed algorithm to achieve an accuracy of $10^{-16}$
when calculating the integral uses no more than 90 points
(\(n \leq 90\)) with \(\omega = 0.1,\ldots 100\). A significant gain in
the number of addition/subtraction and multiplication/division
operations is achieved when the frequency value is greater than the
number n, which ensures the diagonal dominance of the system of linear
algebraic equations in the matrix~(\ref{eq:15}).

It is useful to compare the algorithm we developed for finding the
integrals of highly oscillating functions with the results of~\cite{10}
which presents various and carefully selected numerical examples for
various classes of amplitude functions.

\subsection{Example 3}

Consider the calculation of the integral with an exponential (integer)
function as the amplitude
\begin{equation}
  \label{eq:20}
  I\left( \alpha,\omega \right) =
  \int_{-1}^{1}{e^{\alpha(x - 1)}
    e^{{i\omega x}}{dx}},
\quad
\alpha = 16,64; \omega = 20,1000.
\end{equation}

\begin{figure}
  \centering
  \includegraphics[width=0.8\linewidth]{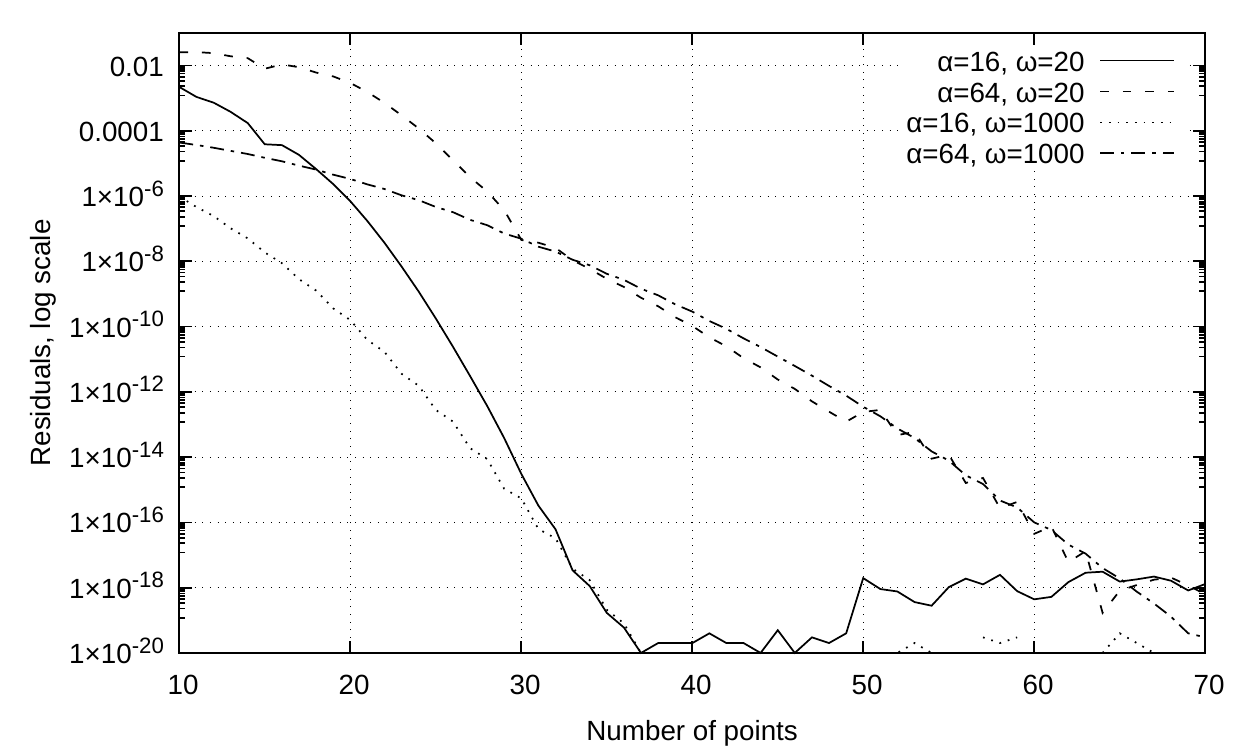}
  \caption{Plot of the absolute error of the approximation of the integral~(\ref{eq:20}) with \(\alpha = 16, 64\); at \(\omega = 20\) and
    \(\omega = 1000\) depending on the number of nodes of the collocation
    method}
  \label{fig:3}
\end{figure}

The exact value of the integral can be calculated by the formula
\(I\left( \alpha,\omega \right) = \frac{2 e^{- \alpha}\sinh(\alpha +
  i\omega)}{(\alpha + i\omega)}\)~\cite{10}.  The plot of the
deviation of the integral calculated by us from the exact one
depending on the number of collocation points (absolute error) is
shown in Fig~\ref{fig:3}.

Comparison with the results of~\cite{10} shows that the accuracy of
calculating the integrals practically coincides with that
of~\cite{10}.  Our advantage is the much simpler form of the matrix of
a system of linear equations. In the best case, when
\(\left| \omega \right| > n\) is a triangular matrix with a main
diagonal and two upper codiagonals.  If \(n > \left| \omega \right|\),
then the transition to the search for a normal solution to a system
with a positively defined Hermitian five-diagonal matrix allows us to
create a numerically stable solution scheme. In both cases, the
LU-decomposition method for SLAEs with band matrices demonstrated the
best results in accuracy and speed.

Below are similar plots, but in the case of calculation using the
preliminary conversion of the triangular matrix of the system to a band
form.

\subsection{Example 4}

In this example~\cite{10}, the highly oscillating function
\(e^{i2\pi\alpha x}\) is considered as the amplitude one. It is clear
that in this case, to achieve the same accuracy in calculating the
integral as in the previous example, a larger number of collocation
points will be required.

\begin{equation}
  \label{eq:21}
  I\left( \alpha,\omega \right) =
  \int_{-1}^{1}{e^{i2\pi\alpha x}e^{{i\omega x}}{dx}},
  \quad
  \alpha = 5,10; \omega = 20,1000.
\end{equation}

\begin{figure}
  \centering
  \includegraphics[width=0.8\linewidth]{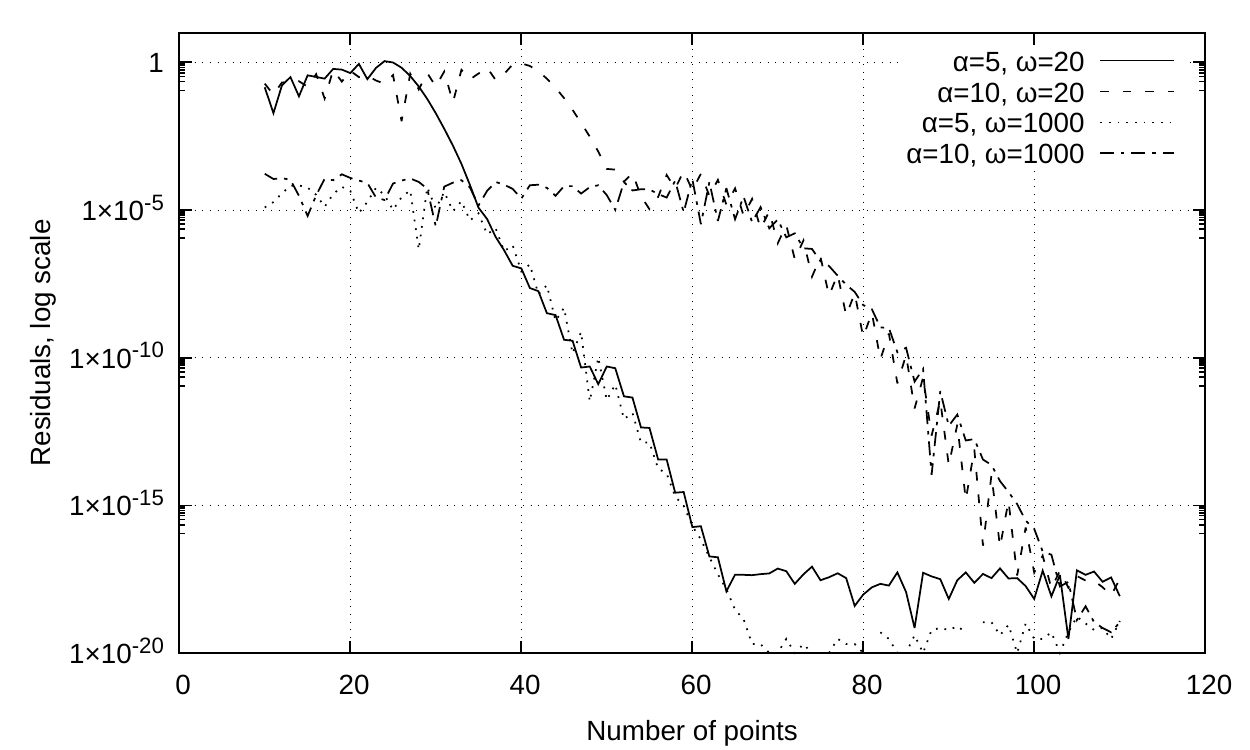}
  \caption{The graph of the absolute error of the approximation of the
    integral~(\ref{eq:21}) with \(\alpha = 5,10\), \(\omega = 20, 1000\) depending on
    the number of nodes of the collocation method}
  \label{fig:4}
\end{figure}

Figure~\ref{fig:4} shows the dependence of the absolute error on the number of
collocation nodes. Similar, in comparison with the previous example,
accuracy (of the order of \(10^{- 17}\)) of the deviation of the
calculated value of the integral from the exact value is achieved only
for $n>100$ in the case \(\alpha = 10\) for various frequency
values \(\omega = 20,\ 1000\).

Calculation of this integral using various methods (triangular
decomposition and reduction to band form) did not reveal a significant
difference in the results. The solutions of the corresponding systems of
linear algebraic equations based on expressions~(\ref{eq:15}) and~(\ref{eq:16})
practically coincide and give the same results.

\subsection{Example 5}

The example demonstrates the calculation of the integral in the case
when the amplitude function is the generating function of the Chebyshev
polynomials of the first kind.

\begin{equation}
  \label{eq:22}
  I\left( \alpha,\omega \right) =
  \int_{-1}^{1}
  \frac{1 - \alpha^{2}}{1 - 2\alpha x + \alpha^{2}
    e^{{i\omega x}}{dx}},
  \quad
  \alpha = 0.8, 0.9; \omega = 20, 1000.
\end{equation}

The behaviour of the amplitude function should lead to an almost linear
dependence of the approximation accuracy on the number of points for
various values of the parameters \(\alpha\) and \(\omega\). Numerical
experiments carried out confirm this assertion.

\begin{figure}
  \centering
  \includegraphics[width=0.8\linewidth]{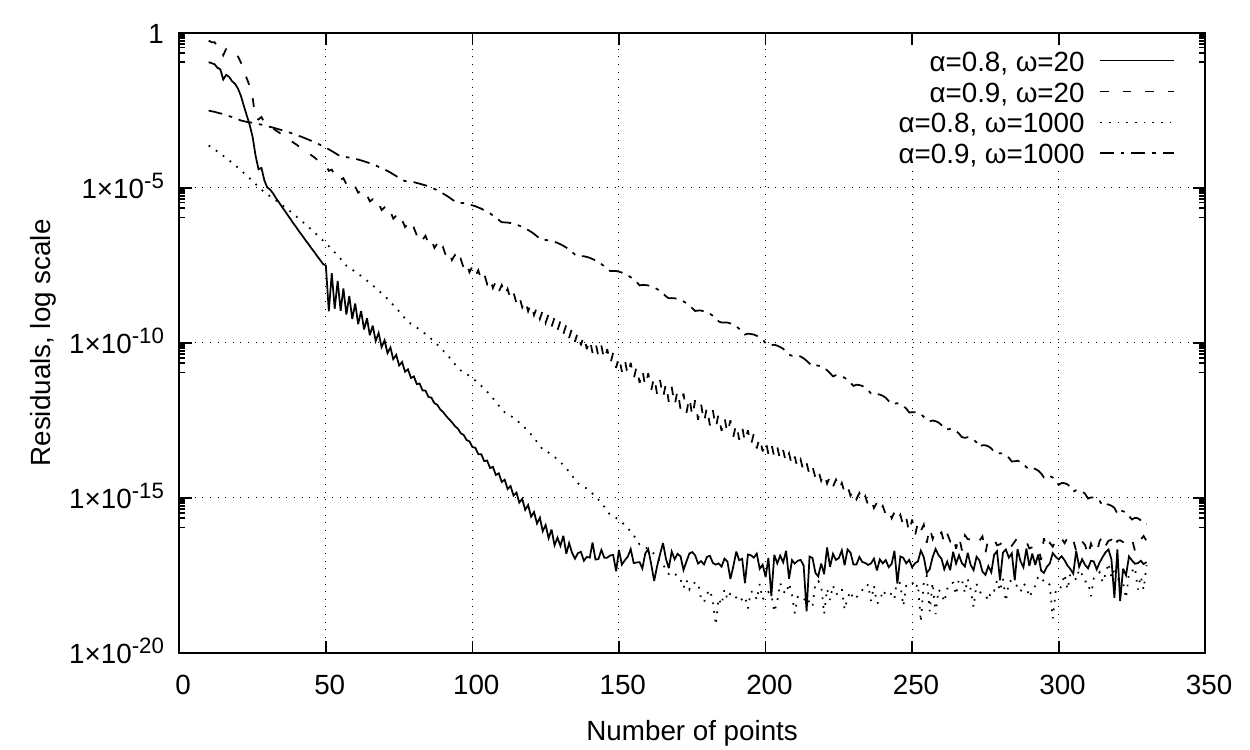}
  \caption{The plot of the absolute error of approximation of the
    integral~(\ref{eq:22}) with
  \(\alpha = 0.8, 0.9\), \(\omega = 20, 1000\) depending on the
  number of nodes of the collocation method}
  \label{fig:5}
\end{figure}

Moreover, the accuracy of calculating the integrals is not inferior to
the accuracy of the methods of \cite{10}

\subsection{Example 6}

Amplitude is a bell-shaped function
\begin{equation}
  \label{eq:23}
  I\left( \alpha,\omega \right) =
  \int_{-1}^{1}{\frac{1}{x^{2} + \alpha^{2}}e^{{i\omega x}}{dx}},
  \quad
  \alpha = 1/4, 1/8; \omega = 20,1000.
\end{equation}

The example is rather complicated for interpolation by Chebyshev
polynomials. To achieve acceptable accuracy (\(10^{- 18}\)), the
deviation of the calculated value of the integral from the exact one
requires about 300 approximation points both for small values of
\(\omega = 20\) and for large \(\omega = 1000\).

\begin{figure}
  \centering
  \includegraphics[width=0.8\linewidth]{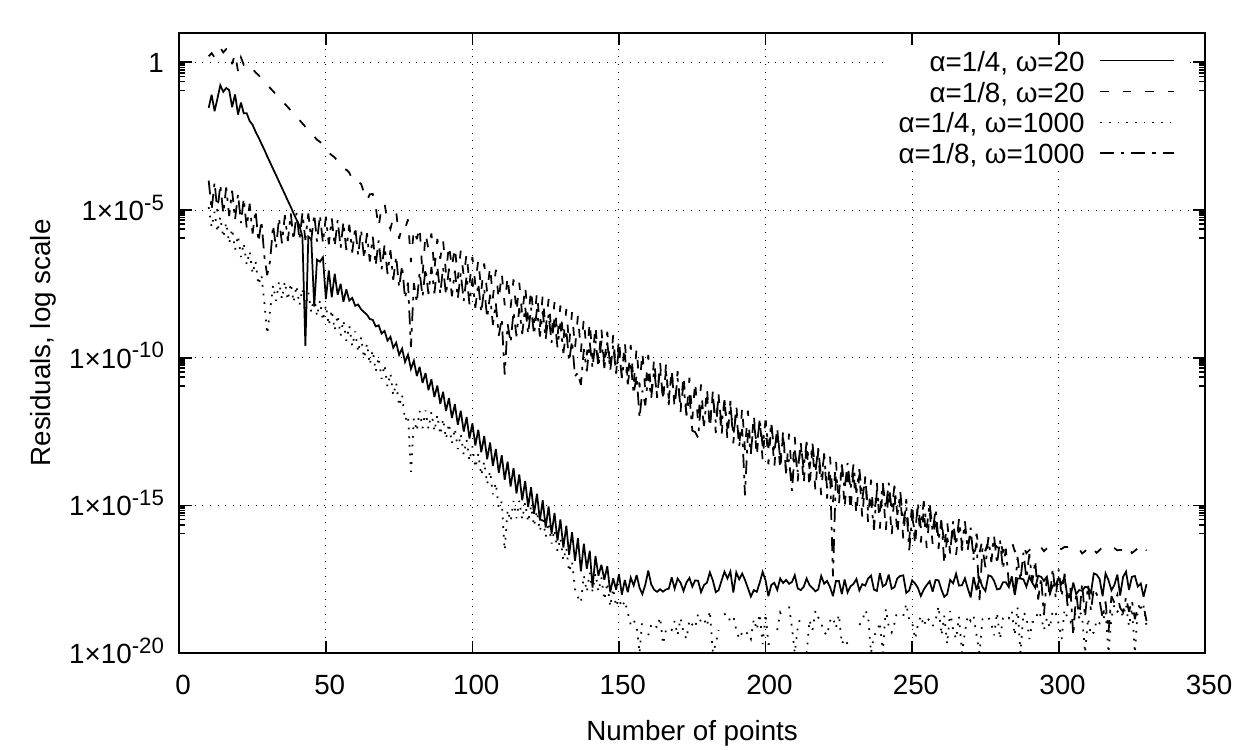}
  \caption{Plot of the absolute error of approximation of the integral~(\ref{eq:23})
    with \(\alpha = 1/4, 1/8\), \(\omega = 20,\ 1000\) depending on the
    number of nodes of the collocation method}
  \label{fig:6}
\end{figure}

Plot of the absolute error of approximation of the integral~(\ref{eq:23})
with \(\alpha = 1/4,\ 1/8\); \(\omega = 20,\ 1000\) depending on the
number of nodes of the collocation method. Logarithmic scale.

\subsection{Example 7}

We give an example of integration when the amplitude function has
second-order singularities at both ends of the integration interval
\begin{equation}
  \label{eq:24}
  I\left( \omega \right) =
  \int_{-1}^{1}{{(1 - t^{2})}^{3/2}e^{{i\omega x}}{dx}},
  \quad \omega = 20,1000.
\end{equation}

The value of this integral can be calculated in an analytical form:
\(I(\omega) = 3\pi J_{2}(\omega)/\omega^{2}\). We present
the numerical values of the integral for various values of the
frequency:
$$
\begin{gathered}
I (20 ) = - 0.00377795409950960,\\
I (1000 ) = - 2.33519886790130 \times 10^{-7}.
\end{gathered}
$$

\begin{figure}
  \centering
  \includegraphics[width=0.8\linewidth]{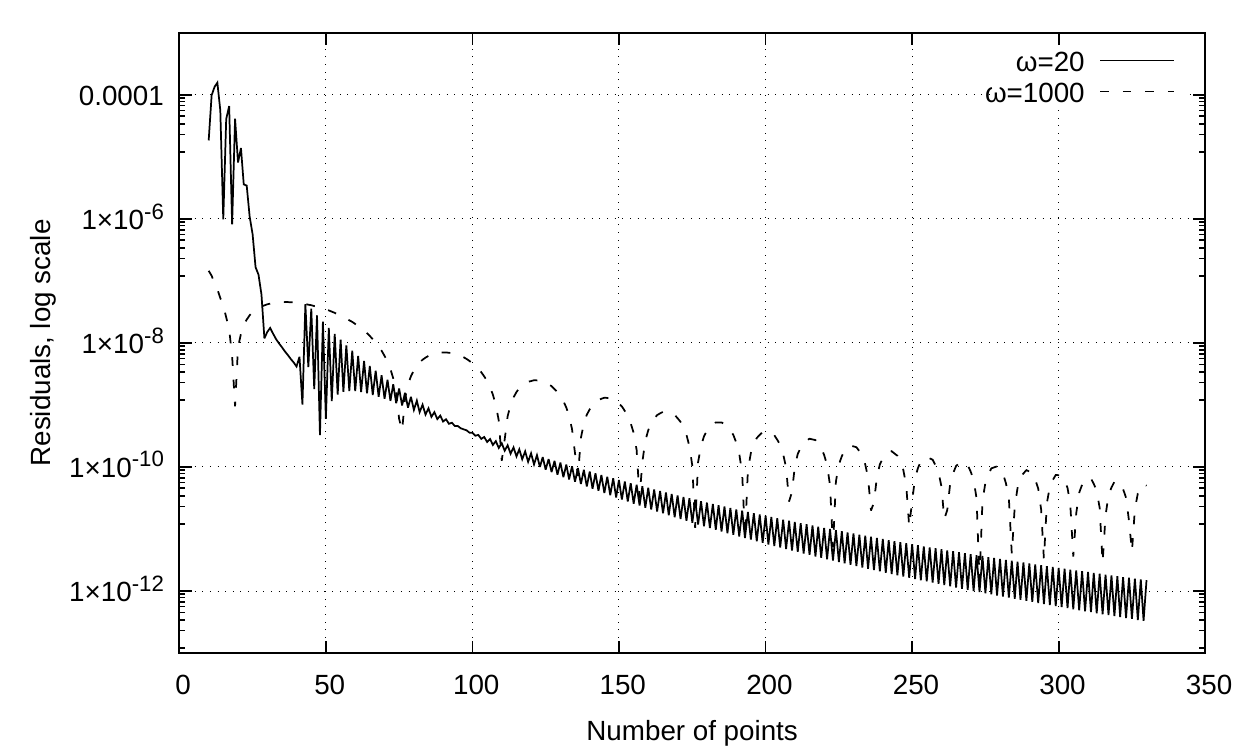}
  \caption{The plot of the absolute error of the approximation of the
    integral~(\ref{eq:24}) with \(\omega = 20,\ 1000\) depending on the number of
    nodes of the collocation method}
  \label{fig:7}
\end{figure}

Similar to the previous example, to achieve good accuracy in calculating
the integral, it is necessary to consider a large number of collocation
points. However, for this type of amplitude functions, the method
presented in the article works reliably both in the case of low and high
frequencies.

The given examples demonstrate that the dependence of the solution on
the number of approximation points is similar to the dependence
demonstrated in Olver \cite{20} and Hasegawa \cite{10}. The advantage of
our approach is the simplicity of the algorithm and the high speed of
solving the resulting very simple system of linear algebraic equations.
If it is necessary to repeatedly integrate various amplitude functions
at a constant frequency, multiple gains are possible due to the use of
the same LU-decomposition backtracking procedure.

\section{Conclusion}
\label{sec:conclusion}

A simple, effective, and stable method for calculating the integrals of
highly oscillating functions with a linear phase is proposed. It is
based on Levin's brilliant idea, which allows the use of the collocation
method to approximate the antiderivative of the desired integral. The
original formulation of the problem by Levin and his followers suggests
a possible ambiguity in finding the antiderivative. Using the expansion
in slowly oscillating polynomials provides a slowly changing solution of
the differential equation.

The transition from a solution in physical space to a solution in
spectral space makes it possible to effectively use the discrete
orthogonality property of the Chebyshev mapping matrix on a
Gauss-Lobatto grid. With this transformation, the uniqueness of the
solution of the studied system is preserved, and its structure from a
computational point of view becomes easier.

There are a large number of works using various approaches aimed to
offer fast and effective methods for solving SLAEs that arise when
implementing the collocation method. However, many methods \cite{8,10}
encounter instability when solving the corresponding systems of linear
equations. When using Chebyshev differentiation matrices in physical
space, instability is explained primarily by the degeneracy of these
matrices and the huge spread of eigenvalues of the matrix of the
collocation method system. The approach to solving the differential
equation based on the representation of the solution, as well as the
phase and amplitude functions, in the form of expansion in finite series
by Chebyshev polynomials and the use of three-term recurrence relations
\cite{8,10,11} also does not provide a stable calculation for
\(n > |\omega|\). To overcome instability, various methods of
regularizing the systems under study are proposed.

In our work, we propose a new method for improving computational
properties by preconditioning (preliminary multiplication by a
non-degenerate band matrix) of the system in the spectral representation
and by searching for its pseudo-normal solution. The proposed method has
been reduced to solving a SLAE with a triangular band matrix or (in the
worth case when \(n > |\omega|\)) Hermitian five-diagonal matrix.
That approach provides a significant reduction in the number of
operations. Several numerical examples demonstrate the advantages of the
proposed effective stable numerical method for integrating highly
oscillating functions with a linear phase.

\begin{acknowledgments}

This paper has been supported by the RUDN University Strategic Academic Leadership Program
(Konstantin P. Lovetskiy, implementation
of the computer code and supporting algorithms).  The reported study
was funded by Russian Foundation for Basic Research (RFBR), project
number 18-07-00567 (Leonid A. Sevastianov, oversight and leadership
responsibility for the research activity planning and execution).
The reported study was funded by Russian Foundation for Basic
Research (RFBR), project number 19-01-00645 (Dmitry S. Kulyabov,
preparation of the published work, specifically visualization and
data presentation).%

\end{acknowledgments}

 \bibliographystyle{elsarticle-num-names}

\bibliography{bib/linear-phase/cite}

\end{document}